\documentclass[reqno,11pt]{amsart}

\usepackage{color}

\newtheorem{theorem}{Theorem}[section] 
\newtheorem{lemma}[theorem]{Lemma}
\newtheorem{corollary}[theorem]{Corollary}
 \theoremstyle{remark}
\newtheorem{remark}[theorem]{Remark}

\newcommand\cD{\mathcal{D}}
\newcommand\cL{\mathcal{L}}

\newcommand\bL{\mathbb{L}}

\newcommand\bR{\mathbb{R}}

\newcommand\Lie{{\rm Lie\,}}

\newcommand{\mysection}[1]{\section{#1}
      \setcounter{equation}{0}}

\theoremstyle{definition}
\newtheorem{assumption}{Assumption}[section]
\newtheorem{definition}{Definition}[section]

\newcommand\cbrk{\text{$]$\kern-.15em$]$}}
\newcommand\opar{\text{\raise.2ex\hbox{${\scriptstyle | }$}\kern-.34em$($} }
\newcommand{\tr}{\text{\rm tr}\,}

 \makeatletter
 \def\dashint{%
 \operatorname%
 {\,\,\text{\bf--}\kern-.98em\DOTSI\intop\ilimits@\!\!}}
 \makeatother

\newcommand\supp{\text{\rm supp}\,}

\newcommand\cB{\mathcal{B}}

\newcommand\cF{\mathcal{F}}

\newcommand\cP{\mathcal{P}}

\newcommand\frD{\mathfrak{D}}

\newcommand\frM{\mathfrak{M}}

\renewcommand\det{{\rm det}\,}

\begin{document}

\title[H\"ormander's theorem]{H\"ormander's theorem for 
stochastic partial differential equations} 
\author{N.V. Krylov}
\thanks{The  author was partially supported by
 NSF Grant DMS-1160569}
\email{krylov@math.umn.edu}
\address{127 Vincent Hall, University of Minnesota,
 Minneapolis, MN, 55455}

\begin{abstract}
We prove H\"ormander's type hypoellipticity theorem for 
stochastic partial differential equations
when the coefficients are only measurable
with respect to the time variable.
The need for such kind of results comes from filtering
theory of partially observable diffusion processes,
when even if the initial system is autonomous,
the observation process enters the coefficients
of the filtering equation and makes them time-dependent
with no good control on the smoothness of the coefficients
with respect to the time variable.
 
\end{abstract}

\keywords{Hypoellipticity, SPDEs, H\"ormander's theorem}

\subjclass[2010]{60H15, 35R60}

\maketitle

\mysection{Introduction}
                                            \label{section 6.28.1}

Let $(\Omega,\cF,P)$ be a complete probability space
with an increasing filtration $\{\cF_{t},t\geq0\}$
of complete with respect to $(\cF,P)$ $\sigma$-fields
$\cF_{t}\subset\cF$. Let $d_{1}\geq1$ be an integer and let
 $w^{k}_{t}$, $k=1,2,...,d_{1}$, be independent one-dimensional
Wiener processes with respect to $\{\cF_{t}\}$.

Fix an integer $d\geq1$ and introduce $\bR^{d}$
as a Euclidean space of column-vectors (written in a common
abuse of notation as) $x=(x^{1},...,x^{d})$.
Denote 
$$
D_{i}=\partial/\partial x^{i},\quad D_{ij}=D_{i}D_{j}
$$
and 
for an $\bR^{d}$-valued function 
$\sigma_{t}(x)=\sigma_{t}(\omega,x)$ on $\Omega\times[0,
\infty)\times\bR^{d}$ 
and functions $u_{t}(x)=u_{t}(\omega,x)$
on $\Omega\times[0,\infty)\times\bR^{d}$ set
$$
L_{\sigma_{t}}u_{t}(x)=[D_{i}u_{t}(x)]\sigma^{i}_{t}(x).
$$

Next take an integer $d_{2}\geq 1$,
assume that we are given $\bR^{d}$-valued functions
$\sigma^{k}_{t}=(\sigma^{ik}_{t})$, $k=0,...,d_{1}+d_{2}$, on 
$\Omega\times[0,\infty)\times\bR^{d}$,
which are infinitely differentiable with respect to $x$
for any $(\omega,t)$, and define the 
operator 
\begin{equation}
                                                       \label{7.22.4}
L_{t}=(1/2)\sum_{k=1}^{d_{2}+d_{1}}L^{2}_{\sigma^{k}_{t}}
+L_{\sigma^{0}_{t}}.
\end{equation}
Assume that on $\Omega\times[0,\infty)\times\bR^{d}$
we are also given certain real-valued 
functions $c_{t}(x)$ and $\nu^{k}_{t}(x)$, $k=1,...,d_{1}$,
which are infinitely differentiable with respect  to $x$,
and that on $\Omega\times[0,\infty)\times
\bR^{d}$ we are given real-valued
functions $f_{t}$ and $g_{t}^{k}$,  $k=1,...,d_{1}$. Then under natural
additional assumptions which will be specified later the
SPDE
\begin{equation}
                                                         \label{6.28.1}
du_{t}=(L_{t}u_{t}+c_{t}u_{t}+f_{t})\,dt
+(L_{\sigma^{ k}_{t}}u_{t}+\nu^{k}_{t}u_{t}+g^{k}_{t})\,dw^{k}_{t}
\end{equation}
makes sense (where and below the summation convention over
repeated indices is enforced regardless of whether they
stand at the same level or at different ones).

The main goal of this paper is to show, somewhat loosely speaking,
that, if $\Omega_{0}\in\cF$,  $ (s_{1},s_{2})\in(0,\infty)$ 
and for any $ \omega\in\Omega_{0}$
and $t\in(s_{1},s_{2})$
the Lie algebra generated by the vector-fields $\sigma^{d_{1}+k}_{t}$,
$k=1,...,d_{2}$, has dimension $d$ everywhere in a ball $B$
in $\bR^{d}$
and $f_{t}$ and $g_{t}^{k}$ are infinitely differentiable
in $B $ for any $ \omega\in\Omega_{0}$
and $t\in(s_{1},s_{2})$, then any function
$u_{t}$ satisfying \eqref{6.28.1} in $\Omega_{0}\times
(s_{1},s_{2})\times B $, for almost any $\omega\in\Omega_{0}$,
coincides on 
$(s_{1},s_{2})\times B $ with a function which is infinitely
differentiable with respect to $x$.
Thus, under a {\em local\/} H\"ormander's type condition
we claim the local hypoellipticity of the equation.

 It is worth
mentioning   article \cite{CM}
where the authors prove the hypoellipticity
for SPDEs whose coefficients do not
explicitly depend on time and $\omega$
under  H\"ormander's type condition which is {\em global}
 but otherwise much weaker    
 than ours.  The  dependence on the time variable $t$  and $\omega$
of the coefficients in \cite{CM} is allowed 
only through an argument in which   a Wiener process
is substituted. However, it seems to the author of the present article
that there is a gap
in the arguments in \cite{CM} when the authors claim that
one can estimate derivatives of order $s+\varepsilon$ ($\varepsilon>0$)
of solutions through derivatives of order $s$
for any $s\in(-\infty,\infty)$ and not only
for $s=0$. The claim is only proved for $s=0$ in \cite{CM}
 and even if there
are no stochastic terms the proof of the claim is not
completely trivial (see the comment below
formula (5.2) in \cite{Horm}).  It is worth noting that our methods
are absolutely different from the methods in \cite{CM}.
Our main method of proving Theorems \ref{theorem 6.30.1}
and \ref{theorem 7.22.1} is based on an observation
by A. Wentzell \cite{We65}
who discovered the It\^o-Wentzell formula and used it to make
a random change of coordinates in such a way that the stochastic terms
in the transformed equation disappear so that
we can use the results from \cite{Horm}.
 We apply this method locally.

Kunita in \cite{Ku81} also uses Wentzell's reduction
of SPDEs with even time-inho\-mogeneous
coefficients to deterministic equations with random
and time-dependent coefficients  
satisfying a {\em global\/}  H\"ormander's type condition.
He writes that the 
probabilistic
approach to proving H\"ormander's
theorem developed by Malliavin \cite{Ma}, Ikeda and Watanabe \cite{IW},
 Stroock \cite{St},
and Bismut \cite{Bi} can be applied to the case of   operators
continuously depending on the time parameter $t$.
In \cite{Ku82} he replaces this list of references with 
\cite{Ma}, \cite{IW}, \cite{St2}, and \cite{Bi2}.
However, to the best of the author's knowledge until now the best result
in proving H\"ormander's theorem {\em by using the Malliavin calculus\/}
for parabolic
equations with the coefficients only continuous 
with respect to $t$ are obtained in
\cite{CaM} where equations with coefficients that are
{\em H\"older} continuous 
in $t$ 
are considered. In our case the coefficients
are only assumed to be predictable, so that if they
are not random, then  their measurability  with respect to $t$
suffices. Another objection against the arguments in 
\cite{Ku81} and \cite{Ku82} is that the reduction of SPDEs is
done globally and yields deterministic parabolic equations
with random coefficients without any control
on their behavior as $|x|\to\infty$, which is needed for
any existing theory of unique solvability of such equations.

Wentzell's method allows us to derive from
 a local version  of H\"ormander's type condition
    infinite differentiability of solutions
at the same locality,
whereas in \cite{CM}, \cite{Ku81}, and \cite{Ku82}
a global condition is imposed and the way $\omega$ and $t$
enter the coefficients
is quite restrictive.
Another difference between our results and those
in \cite{CM} is that we prove infinite differentiability of
any generalized solution and not only of measure-valued ones.

Speaking about generalized solution,
our functions $u_{t},f_{t},g^{k}_{t}$
are, actually, assumed to be given on a subset of $\Omega\times
[0,\infty) $ and take values in $\cD$, which is the space 
of
generalized functions on  $\bR^{d}$. 

One more issue worth noting is that we derive a priori
estimates which will allow us in a subsequent article
not only show that the filtering density for $t>0$
 is in $C^{\infty}$ if
the unobservable process starts at any fixed point $x$ but
also prove that it is infinitely differentiable
with respect to $x$. As far as the author is aware
such kind of results was never proved for degenerate SPDEs.

We finish the introduction with a few more notation
and a description of the structure of the article.
For (generalized) 
functions $u$ on $\bR^{d}$ by $Du$ we mean the {\em row-vector\/}
$(D_{1}u,...,D_{d}u)$ and when we write $Du \phi$ we always mean
$(Du)\phi$. In this notation
$$
L_{\sigma_{t}}u_{t} =[D_{i}u_{t} ]\sigma^{i}_{t} 
= Du_{t}  \sigma_{t} .
$$
One knows that the product of any generalized function
and an infinitely differentiable one is again a generalized function
and that any generalized function is infinitely differentiable
in the generalized sense, so that what is said   above 
has perfect sense.

 For $R, t\in(0,\infty)$   set
$$
B_{R}=\{x\in\bR^{d}:|x|<R\}, 
\quad C_{t,R}=(0,t)\times B_{R},
$$
and denote by $\cD_{R}$ the set of generalized functions
on $B_{R}$.
In the whole article $T,R_{0}$ are fixed numbers
from $ (0,\infty)$.

The rest of the article is organized as follows.
In Section \ref{section 7.16.1} we state our main results,
Theorems \ref{theorem 6.30.1} and \ref{theorem 7.22.1}.
Section \ref{section 9.8.1} contains a computation
of the determinant of a matrix-valued process satisfying
a linear stochastic equation. A very short
 Section \ref{section 7.13.1} reminds the reader
one of properties of stochastic integrals of
Hilbert-space valued processes. In Section \ref{section 6.30.1}
we discuss some facts  related to stochastic
flows of diffeomorphisms and change of variables.
The reader can find  in \cite{Ku90}
much more information about
stochastic
flows of diffeomorphisms in a much more general setting.
 Our discussion is more elementary than in \cite{Ku90}
 albeit it is only valid in a particular case
we need. 
 In Section \ref{section 6.30.4} we prove a version of 
the It\^o-Wentzell formula 
we need. Finally, in
  Sections \ref{section 6.30.2} and \ref{section 8.5.1}
we prove Theorems \ref{theorem 6.30.1} and \ref{theorem 7.22.1},
respectively.

\mysection{Main results}
                                                   \label{section 7.16.1}
Denote by $\cP$ the predictable
$\sigma$-field in $\Omega\times(0,\infty)$
associated with $\{\cF_{t}\}$. 

\begin{definition} 
Denote by $\frD(C_{T,R_{0}})$                   \label{def 10.25.1}
 the set of all $\cD_{R_{0}}$-valued 
functions $u$ (written 
as $u_{t}(x)$ in a common abuse of notation)
on $\Omega\times[0,T]$ such that, for any 
$\phi\in C_{0}^{\infty}(B_{R_{0}})$,
 the restriction of the function $(u_{t},\phi)$ 
on $\Omega\times(0,T]$ is $\cP$-measurable 
and $(u_{0},\phi)$ is $\cF_{0}$-measurable.
For $p=1,2$  
denote by $\mathfrak{D}^{-\infty}_{p }(C_{T,R_{0}})$ 
the subset of $\frD(C_{T,R_{0}})$
consisting of $u$ such that 
for any $\zeta\in C^{\infty}_{0}(B_{R_{0}})$
  there exists an $m\in \bR$ such that for any 
  $\omega\in\Omega$,
 for almost all $t\in[0,T]$,
we have $\zeta u_{t}\in H^{m}_{2} $ ($=(1-\Delta)^{-m/2}
\cL_{2},\,\,\cL_{2}=\cL_{2}(\bR^{d})
 $)
and  
\begin{equation}
                                            \label{11.16.2}
\int_{0}^{T} \| u_{t}  
\zeta \|_{H^{m}_{2} }^{p}\,dt<\infty.
\end{equation}

\end{definition}

\begin{definition} 
                                           \label{def 10.25.3}
Assume that we are given some $u,f,g^{k}
\in\mathfrak{D}(C_{T,R_{0}})$, 
$k=1,...,d_{1}$ (not necessarily those from Section \ref{section 6.28.1}).
 We say that the equality
\begin{equation}
                                           \label{11.16.3}
du_{t}(x)=f_{t}( x)\,dt+
g_{t}^{k}( x)\,dw^{k}_{t},\quad (t,x)\in  C_{T,R_{0}} ,
\end{equation}
holds {\em in the sense of distributions\/} if
$ f \in\mathfrak{D}^ {-\infty}_{1 }(C_{T,R_{0}})$, 
$g^{k}\in
\mathfrak{D}^ {-\infty}_{2 }(C_{T,R_{0}})$, $k=1,...,d_{1}$, and
for
 any $\phi\in C_{0}^{\infty}(B_{R_{0}})$, 
 with probability one  we have 
\begin{equation}
                                             \label{12.23.40}
(u_{t } ,\phi)=(u_{0} ,\phi)+\int_{0}^{t}
(f_{s},\phi)\,ds+\sum_{k=1}^{d_{1}}
\int_{0}^{t}  (g^{k}_{s},\phi)\,dw^{k}_{s}
\end{equation}
for all $t\in [0,T]$, where, as usual, $(\cdot,\cdot)$
stands for pairing of generalized  and test functions.
\end{definition}

\begin{remark}
                                          \label{remark 6.28.2}
Observe that if $g^{k}\in\frD^{-\infty}_{2 }(C_{T,R_{0}})$, $\phi,
\zeta\in C^{\infty}_{0}(B_{R_{0}})$, and $\zeta=1$ on the support of $\phi$,
then
$$
|(g^{k}_{s},\phi)|^{2}=|(\zeta g^{k}_{s},\phi)|^{2}
\leq \|\zeta g^{k}_{s}\|_{H^{m}_{2}}^{2}\|\phi\|_{H^{-m}_{2}}^{2}
$$ 
and the right-hand side has finite integral over $[0,T]$ (a.s.)
if $m$ is chosen appropriately. This and a similar estimate
concerning $(f_{s},\phi)$ shows that the right-hand side in
\eqref{12.23.40} makes sense.

\end{remark}

In the following assumption we are talking about the objects
from Section~\ref{section 6.28.1}.

\begin{assumption}
                                        \label{assumption 6.28.1}
(i) The  functions $\sigma^{k}_{t}(x)$, $k=0,...,d_{1}+d_{2}$,
$c_{t}$, $\nu^{k}_{t}$, $k=1,...,d_{1}$, are infinitely
differentiable with respect to $x$ and each of their derivatives
of any order
is bounded on $\Omega\times[0,T]\times B_{R_{0}}$.
These functions are predictable with respect to $(\omega,t)$
for any $x\in B_{R_{0}}$;

(ii) We have $u,f,g^{k}\in\frD^{-\infty}_{2 }(C_{T,R_{0}})$, $k=1,...,d_{1}$;

(iii) Equation \eqref{6.28.1} holds on 
$ C_{T,R_{0}} $ in the sense of
Definition \ref{def 10.25.3}.

(iv) for any $\zeta\in C^{\infty}_{0}(B_{R_{0}})$
  there exists an $m\in \bR$ such that for any 
  $\omega\in\Omega$, we have $u_{0}\zeta\in H^{m}_{2}$.

\end{assumption}

\begin{remark}
                                        \label{remark 11.1.1}
The argument in Remark \ref{remark 6.28.2}
shows that \eqref{6.28.1} has perfect sense owing to
Assumptions \ref{assumption 6.28.1} (i), (ii), and we need
$u \in\frD^{-\infty}_{2 }(C_{T,R_{0}})$ in contrast with Definition
\ref{def 10.25.3} because $Du$ and $u$ enter the stochastic
part in \eqref{6.28.1}.  

Furthermore, under Assumption \ref{assumption 6.28.1}
for any $\zeta\in C^{\infty}_{0}(B_{R_{0}})$ there is an
$m$ such that $u_{0}\zeta\in H^{m}_{2}$ and
$$
\int_{0}^{T} \big(\| u_{t}  
\zeta \|_{H^{m}_{2} }^{2}+\| f_{t}  
\zeta \|_{H^{m}_{2} }^{2}
+\sum_{k=1}^{d_{1}}\| g^{k}_{t}  
\zeta \|_{H^{m}_{2} }^{2}\big)\,dt<\infty.
$$
It follows by a classical continuity result that
(a.s.) $u_{t}  
\zeta$ is a continuous $H^{m-1}$-valued function on $[0,T]$.
If we drop Assumption \ref{theorem 7.22.1} (iv), then the
same will be true with $(0,T]$ in place of $[0,T]$
since $u_{t}  
\zeta \in H^{m}_{2}$ for almost all $t\in(0,T)$.

\end{remark}

Next, as usual, for two smooth $\bR^{d}$-valued functions $\sigma,\gamma$
on $\bR^{d}$ we set
$$
[\sigma,\gamma]=D\gamma \sigma-D\sigma\gamma,
$$
where for instance $D\gamma$ is the matrix with entries $
(D\gamma)^{ij}=D_{j}\gamma^{i}$, so that
$$
[\sigma,\gamma]^{i}=\sigma^{j}D_{j}\gamma^{i} -\gamma^{j}D_{j}\sigma^{i}.
$$

Then introduce collections of $\bR^{d}$-valued functions
defined on $\Omega\times[0,T]\times B_{R_{0}}$ inductively as
$\bL_{0}=\{\sigma^{d_{1}+1},...,\sigma^{{d_{1}+d_{2}}}\}$,  
$$
\bL_{n+1}=\bL_{n}\cup\{[\sigma^{d_{1}+k},M]:k=1,...,d_{2},
M\in\bL_{n}\},\quad n\geq0.
$$ 

For any multi-index $\alpha=(\alpha_{1},...,\alpha_{d})$,
$\alpha_{i}\in\{0,1,...\}$, introduce as usual
$$
D^{\alpha}=D^{\alpha_{1}}_{1}\cdot...\cdot D^{\alpha_{d}}_{d},
\quad |\alpha|=\alpha_{1}+...+\alpha_{d}.
$$
Also define
 $BC^{\infty}_{b}$ as the set of real-valued  
measurable
 functions $a$ on $\Omega\times[0,T]\times\bR^{d}$
such that, for each $t\in[0,T]$ and $\omega\in\Omega$,
 $a_{t}(x)$ is infinitely differentiable with respect to $x$, and
  for any $\omega\in\Omega$ and multi-index $\alpha$ we have
$$
\sup_{ (t,x)\in [0,T]\times\bR^{d}}|D^{\alpha}a_{t}( x)|<\infty.
$$

  Finally we denote by
$\Lie_{n}$ the set of (finite) linear combinations
of elements of $\bL_{n}$ with  
coefficients which  are
of class $BC^{\infty}_{b}$.
Observe that the vector-field $\sigma^{0}$ is {\em not\/}
explicitly included
into $\Lie_{n}$. Finally, fix $\Omega_{0}\in\cF$, $S \in[0,T)$,
and introduce
 $$
G=(S ,T)\times B_{R_{0}}.
$$

\begin{assumption}
                                             \label{assumption 1.5.1}
For every $\omega\in\Omega_{0}$, 
   $\eta\in C^{\infty}_{0}(S ,T)$,   
  and $\zeta\in C^{\infty}_{0}(B_{R_{0}})$
there exists an $n\in\{0,1,...\}$ such that
 we have  $\xi\eta \zeta  \in\Lie_{n}$ for any $\xi\in\bR^{d}$. 
\end{assumption}

Here is our first main result which is proved in
Section \ref{section 6.30.2}.
We remind the reader that the common way of saying that
a generalized function in a domain is smooth
means that there is a smooth function which, 
as
a   generalized function, coincides
with the given generalized one in the domain
under consideration.

\begin{theorem}
                                               \label{theorem 6.30.1}
Assume that for any $\omega\in\Omega_{0}$, $n=1,2,...$, and $\zeta\in C^{\infty}_{0}(G)$,
for almost any $t\in[S ,T]$ we have $f_{t} \zeta  
\in H^{n}_{2}$  and
$$
\int_{S }^{T} \|f_{t}\zeta \|_{H^{n}_{2}}^{2}
 \,dt<\infty 
$$
and for any $\omega\in\Omega $, $n=1,2,...$, and $\zeta\in C^{\infty}_{0}(G)$,
for almost any $t\in[S ,T]$ we have $g^{k}_{t}\zeta  
\in H^{n}_{2}$, $k=1,...,d_{1}$, and
$$
\sum_{k=1}^{d_{1}}\int_{S }^{T} 
\|g^{k}_{t}\zeta \|_{H^{n}_{2}}^{2} 
\,dt<\infty.
$$
  
Then, for almost all $\omega\in\Omega_{0}$,
  $u_{t}( x)$ is infinitely differentiable 
with respect to $x$ for
$(t,x)\in G$ and each derivative is a continuous function
in $G$.

Furthermore, let $[s_{0},t_{0}]\subset(S ,T)$,
 $r\in(0,R_{0})$,   take a $\zeta\in C^{\infty}_{0}(G)$
such that $\zeta=1$ on a neighborhood of
 $[s_{0},t_{0}]\times \bar B_{r}$, and take
an  $m$ (which exists by definition) such that
\eqref{11.16.2} holds with $p=2$. Then,
for any multi-index $\alpha$ and $l$ such that
\begin{equation}
                                                 \label{7.16.1} 
2(l-|\alpha|-2)>d+1 
\end{equation}
there exists a (random, finite)
constant $N$, independent of $u,f$, and $g^{k}$, such that,
for almost any $\omega\in\Omega_{0}$,
\begin{equation}
                                                 \label{7.22.2}
\sup_{(t,x)\in[s_{0},t_{0}]\times B_{r}}
|D^{\alpha}u_{t}(x)|^{2}\leq N
\int_{S }^{T}\big[\|f_{t}\zeta \|_{H^{l}_{2}}^{2}
+\| u_{t}  
\zeta \|_{H^{m}_{2} }^{2}\big]\,dt,
\end{equation}
provided that $g^{k}_{t}\zeta I_{\Omega_{0}}\equiv0$, $k=1,...,d_{1}$.

\end{theorem}

Here is a result which is ``global'' in $t$.
We derive it from Theorem \ref{theorem 6.30.1}
in Section \ref{section 8.5.1}.

\begin{theorem}
                                 \label{theorem 7.22.1} 
Suppose that 
a stronger assumption than
Assumption \ref{assumption 1.5.1}
is satisfied: For every $\omega\in\Omega_{0}$  
  and
 $\zeta\in C^{\infty}_{0}(B_{R_{0}})$
there exists an $n\in\{0,1,...\}$ such that
 we have  $\xi I_{[S ,T]}\zeta  
\in\Lie_{n}$ for any $\xi\in\bR^{d}$.
Also suppose that the assumption stated in Theorem \ref{theorem 6.30.1}
is satisfied for with
$\zeta\in C^{\infty}_{0}(B_{R_{0}})$
 rather than
$\zeta\in C^{\infty}_{0}(G)$. 

 Then the first assertion of Theorem 
\ref{theorem 6.30.1} holds true with $(S ,T]\times B_{R_{0}}$
in place of $G$, and the second assertion holds with
$s_{0}\in(S ,T)$,
 $t_{0}=T$, and $\zeta\in C^{\infty}_{0}(B_{R_{0}})$,
which equals one in a neighborhood of $\bar B_{r}$.

If we additionally assume that $u_{S } $ is infinitely differentiable
in $B_{R_{0}}$ for every $\omega\in\Omega_{0}$,
then the first assertion of Theorem 
\ref{theorem 6.30.1} holds true with $[S ,T]\times B_{R_{0}}$
in place of $G$, and the second assertion holds with
$s_{0}=S $, $t_{0}=T$, and $\zeta\in C^{\infty}_{0}(B_{R_{0}})$,
which equals one in a neighborhood of $\bar B_{r}$,
 if we add to the right-hand side
of \eqref{7.22.2} a constant (independent of $u$) times
$\|\zeta u_{S }\|^{2}_{H^{l+1}_{2}}$.

\end{theorem}

\begin{remark}
The reader will see that Assumption \ref{assumption 6.28.1} (iv)
will be used only in the proof of the second assertion
of Theorem \ref{theorem 7.22.1} for $S=0$.

\end{remark}

\mysection{On linear stochastic equations}
                                                   \label{section 9.8.1}

 Let $z_{t}$ be a $d\times d$ matrix-valued
continuous $\cF_{t}$-adapted process
satisfying 
$$
z_{t}=I+\int_{0}^{t}
\alpha^{k}_{s}z_{s}\,dw^{k}_{s}
+\int_{0}^{t}
\beta_{s}z_{s}\,ds,\quad s\geq0,
$$
where $\alpha^{k}_{s}$, $k=1,...,d_{1}$,
and $\beta_{s}$ are bounded predictable
$d\times d$ matrix-valued processes and $I$
is the identity $d\times d$ matrix.
The goal of this section is to prove
the following result, which is probably
known, but the author could not find
an appropriate reference. In any case the proof is short.

\begin{lemma}
                          \label{lemma 6.21.1}
For $s\geq0$
\begin{equation}
                                    \label{6.22.1}
\det z_{t}=\exp\big(\int_{0}^{t}
\tr\alpha^{k}_{s}\,dw^{k}_{s}
+\int_{0}^{t}\big[\tr\beta_{s}-
(1/2)\sum_{k=1}^{d_{1}}
\tr((\alpha^{k}_{s})^{2})\big]\,ds\big).
\end{equation}
\end{lemma}

Proof. Take a $d\times d$ nondegenerate
matrix $A=(A^{ij})$
and consider it as a function of its entries
$A^{ij},i,j=1,...,d$. Then
$\det A$ is also a function of $A^{ij}$. 
One knows that (we write $f_{x}$ to denote the derivative of $f$
with respect to $x$)
$$
(\det A)_{A^{ij}}=B^{ji}\det A,
$$
where $B=A^{-1}$. Also as with derivatives
with respect to any parameter
$$
B_{A^{rp}}=-BA_{A^{rp}}B.
$$
Observe that $A_{A^{rp}}^{nm}
=\delta^{rn}\delta^{pm}$. It follows that
$$
B_{A^{rp}}^{ji}=-B^{jn}\delta^{rn}\delta^{pm}
B^{mi}=-B^{jr}B^{pi},
$$
$$
(\det A)_{A^{ij}A^{rp}}=-B^{jr}B^{pi}\det A
+B^{ji}B^{pr}\det A.
$$

Now we can use It\^o's formula. Denote
$x_{t}=z^{-1}_{t}$. Then
$$
d\det z_{t}=x_{t}^{ji}\alpha^{ink}_{t}z^{nj}_{t}
\det z_{t}\,dw^{k}_{t}+x_{t}^{ji}
\beta^{in}_{t}z^{nj}_{t}
\det z_{t}\,dt
$$
$$
+(1/2)\big[x^{ji}_{t}x^{pr}_{t}-x^{jr}_{t}
x^{pi}_{t}]
\alpha^{ink}_{t}z^{nj}_{t}
\alpha^{rmk}_{t}z^{mp}_{t}\det z_{t}\,dt.
$$
We note that
$$
x_{t}^{ji} z^{nj}_{t}=\delta^{in},
\quad  x^{pr}_{t}z^{mp}_{t}=\delta^{rm},
\quad x^{jr}_{t}z^{nj}_{t}=\delta^{rn},
\quad x^{pi}_{t}z^{mp}_{t}=\delta^{im}
$$ 
and conclude that
$$
d\det z_{t}=\det z_{t}\big[\tr \alpha^{k}_{t}
\,dw^{k}_{t}+\tr\beta_{t}\,dt
+(1/2)\big(\sum_{k=1}^{d_{1}}(\tr \alpha_{t}^{k})^{2}
-\tr((\alpha^{k}_{t})^{2})\big)\,dt\big].
$$
We see that $\det z_{t}$ satisfies a linear
equation as long as it stays strictly positive.
A unique solution of this equation
which equals one at $t=0$ is given
by the right-hand side of \eqref{6.22.1}, which 
does not vanish for $t\geq0$.  This shows that
\eqref{6.22.1} holds for all $t\geq0$
and the lemma is proved.

\mysection{On stochastic integrals of
Hilbert-space valued processes}
                                           \label{section 7.13.1}
Let $H$ be a separable Hilbert space (in our applications
$H$ is one of $H^{ -n}_{2}$ with large $n>0$). Take a 
predictable $H$-valued process $h_{t},t\in[0,T]$, such that (a.s.)
$$
\int_{0}^{T}\|h_{t}\|_{H}^{2}\,dt<\infty
$$
for any $\omega$
and set $w_{t}=w^{1}_{t}$.

\begin{lemma}
                                                 \label{lemma 7.13.1}
The stochastic integral
$$
\int_{0}^{t}h_{s}\,dw_{s}
$$
has  a (continuous) modification such that, if 
there is a $\phi\in H$, $(s_{0},t_{0})\subset(0,T)$, and $\omega\in\Omega$
for which $(\phi,h_{r}(\omega))_{H}=0$ for $r\in(s_{0},t_{0})$, then
$$
\big(\phi,\int_{0}^{t}h_{s}\,dw_{s}\big)_{H}
$$
is constant on that $\omega$ for $r\in[s_{0},t_{0}]$.
\end{lemma}

The proof of this lemma is achieved immediately
after one recalls that there exists a sequence $n_{k}\to\infty$
and a  $c\in(0,1)$ such that (a.s.) uniformly
on $[0,T]$
$$
\int_{0}^{t}h_{\kappa(n_{k},s+c)-c}\,dw_{s}:=
$$
$$
\sum_{m=1}^{\infty}I_{s\leq t}h_{t_{mk}-c}I_{t_{mk}\leq s+c<
t_{m+1,k}}(w_{t_{m+1,k}-c}-w_{t_{mk}-c})
\to \int_{0}^{t}h_{s}\,dw_{s},
$$
in $H$,
where $t_{mk}=m2^{-n_{k}} $, $\kappa(n,s)=2^{-n}[2^{n}s]$,
and $h_{t}$ is extended as zero outside $[0,T]$.

\mysection{On some random mappings}
                                              \label{section 6.30.1}

Here we suppose that Assumption \ref{assumption 6.28.1} (i)  
is satisfied with $R_{0}=\infty$ and, moreover,
 there is an $R\in(0,\infty)$
such that, for any $k=0,1,...,d_{1}$ and $\omega,t$, we have
$ \sigma^{ k}_{t}(x)=0$
 if $|x|\geq R$.

Consider the equation
\begin{equation}
                                             \label{6.28.4}
x_{t}=x-\int_{0}^{t}\sigma^{k}_{s}( x_{s})\,dw^{k}_{s}
-\int_{0}^{t}b_{t}(x_{s})\,ds,
\end{equation}
where
$$
b_{t}(x)=\sigma^{0}_{t}(x)-(1/2)\sum_{k=1}^{d_{1}}D\sigma^{k}
_{t}(x)\sigma^{k}_{t}(x).
$$
As  follows from \cite{BF} (see \cite{Ku90} for more
advanced treatment of the subject), there exists
a function $X_{t}(x)$ on $\Omega\times[0,T]\times\bR^{d}$ 
such that

(i) it is continuous in $(t,x)$ for any $\omega$
along with each derivative of $X_{t}(x)$ of any order
with respect to $x$,

(ii) it is  $\cF_{t}$-adapted for any $(t,x)$,

(iii) for each $x$ with probability one it satisfies
\eqref{6.28.4} for all $t\in[0,T]$,

(iv) the matrix $DX_{t}(x)$ for any $x$
with probability one satisfies
$$
DX_{t}(x)=I-\int_{0}^{t}D\sigma^{ k}_{s}
(X_{s}(x))DX_{s}(x)\,dw^{k}_{s}
-\int_{0}^{t}Db _{s}
(X_{s}(x))DX_{s}(x)\,ds
$$
for all $t\in[0,T]$.
 
By Lemma \ref{lemma 6.21.1}
we obtain that for any $x$
with probability one
$$
\det DX_{t}(x)=\exp\big(-\int_{0}^{t}\tr
D\sigma^{ k}_{s}
(X_{s}(x))\,dw^{k}_{s}
$$
$$
-\int_{0}^{t}\big[\tr
Db _{s}-(1/2)\sum_{k=1}^{d_{1}}\tr((D\sigma^{ k}_{s})^{2})]
(X_{s}(x))\,ds\big)
$$
for all $t\in[0,T]$. By  formally considering
the system consisting of equation \eqref{6.28.4}
and the ``equation''
$$
y_{t}=y+\int_{0}^{t}\tr
D\sigma^{ k}_{s}
(x_{s})\,dw^{k}_{s}
$$ 
and applying what
is said above, we see that
there exists a function $I_{t}(x)=I_{t}(\omega,x)$ which is
continuous with respect to $(t,x)\in[0,T]\times\bR^{d}$
for each $\omega$ and such that for each $x$ 
$$
I_{t}(x)=\int_{0}^{t}\tr
D\sigma^{ k}_{s}
(X_{s}(x))\,dw^{k}_{s}
$$
with probability one for all $t\in[0,T]$.
Then for each $(t,x)$ 
with probability one
$$
\det DX_{t}(x)=\exp\big(-I_{t}(x)
$$
$$
-\int_{0}^{t}\big[\tr
Db_{s}-(1/2)\sum_{k=1}^{d_{1}}\tr((D\sigma^{ k}_{s})^{2})]
(X_{s}(x))\,ds\big)
$$
and since both parts of these equality are continuous
with respect to $(t,x)$ the equality holds for all
$(t,x)$ at once with probability one.

It follows that, perhaps after modifying $X_{t}(x)$
on a set of probability zero, we may assume that $
\det DX_{t}(x)>0$ for all $(\omega,t,x)$. Also observe
that obviously $X_{t}(x)=x$ for $|x|\geq R$
and $|X_{t}(x)|\leq R$ for $|x|\leq R$. Hence,
there is a  random variable $\varepsilon=\varepsilon(\omega)>0$
such that $
\det DX_{t} \geq \varepsilon$ and
$$
\det[(DX_{t})^{*} DX_{t}] \geq \varepsilon  
$$
for all $(\omega,t,x)$. Since $DX_{t}(x)$ is a bounded function
of $(t,x)$ for each $\omega$, it follows that the smallest
eigenvalue of the symmetric matrix $ (DX_{t})^{*}DX_{t}$
  is bounded below
by a $\delta=\delta(\omega)>0$, that is
\begin{equation}
                                                        \label{7.1.1}
| DX_{t} \xi|^{2} \geq\delta|\xi|^{2}
\end{equation}
for all $(\omega,t,x)$ and $\xi\in\bR^{d}$. 

Now we need the following consequence of \eqref{7.1.1},
which is proved in a much more general case 
of quasi-isometric mappings of Banach spaces
in 
Corollary of Theorem II of \cite{Jo1}
(see also \cite{Jo2}).  
\begin{lemma}
                                         \label{lemma 7.1.1}
For all $(\omega,t)$, the mapping $X_{t}(x)$ of $\bR^{d}$
is one-to-one and onto $\bR^{d}$.
\end{lemma}

Kunita \cite{Ku82} gives a different proof of
Lemma \ref{lemma 7.1.1} in a much more general case
 based on the fact that the 
mapping $X_{t}(x)$ is obviously homotopic to the identity
mapping (but still in his case an additional effort is applied
because $\bR^{d}$ is not compact). Yet another proof provides the following 
result, in which the nondegeneracy of the Jacobian
is not required and which may have an independent interest.

\begin{lemma}
                                         \label{lemma 9.15.1}
Let $D$ be a  connected bounded domain
in $\bR^{d}$ and $X:\bar D\to\bar D$ be a 
continuous mapping which has bounded and continuous
first-order derivatives in $D$. Assume that $X(x)=x$
if $x\in\partial D$, $\det DX(x) $ does not change
sign in $D$, and for any $x_{0}\in  D$ the mapping $X(x)$
is a homeomorphism if restricted to a neighborhood
of $x_{0}$ (for instance, $\det DX(x)>0$ on $D$).
Then the mapping $X$ is one-to-one and onto $\bar D$.
\end{lemma}

Proof. The fact that the mapping is onto
is an easy consequence of the fact that $\partial X(D)
=X(\partial D)=\partial D$.

To prove that $X$ is one-to-one,
for $n=1,2,...$, $i=(i_{1},...,i_{d})$,
$i_{k}=0,\pm1,...$, introduce
$$
C_{i,n}=(i_{1}/2^{n},(i_{1}+1)/2^{n}]
\times...\times (i_{d}/2^{n},(i_{d}+1)/2^{n}].
$$
Take a domain $D'\subset \bar{D}'\subset D$ and observe that,
because of our assumption that $X$ is a local homeomorphism,
 there exists an $n$ such that $X$ restricted
to $C_{i,n}\cap D'$ is one-to-one whenever this intersection
is nonempty. In that case also
$$
\text{\rm Vol\,}X(C_{i,n}\cap D')=\int_{C_{i,n}\cap D'}
|\det DX(x)|\,dx.
$$
Summing  up these relations and then letting $D'\uparrow D$
 we obtain
\begin{equation}
                                             \label{9.15.1}
\text{\rm Vol\,}D=\text{\rm Vol\,}X(D)\leq \int_{ D}
|\det DX(x)|\,dx.
\end{equation}

Note that \eqref{9.15.1} holds without the assumption
that $X$ does not move the points on the boundary of $D$.
Also note for the future that $X$ is Lipschitz continuous
in $\bar{D}$. Indeed, if $x_{1},x_{2}\in\bar{D}$
and the open straight segment connecting $x_{1}$ and $x_{2}$
belongs to $D$, then $|X(x_{1})-X(x_{2})|\leq N_{0}|x_{1}-x_{2}|$,
where $N_{0}$ is the supremum of $\|DX\|$ over $D$. If not the whole
of the segment is in $D$, then
denote by $y_{1}\in\partial D$ and $y_{2}\in\partial D$ the closest
points to $x_{1}$ and $x_{2}$, respectively, 
on the closure of this segment  .
Then 
$$
|X(x_{1})-X(x_{2})|\leq N_{0}|x_{1}-y_{1}|+|y_{1}-y_{2}|+
N_{0}|y_{2}-x_{2}|\leq (N_{0}+1)|x_{1}-x_{2}|.
$$

Next, concentrate on the case that $\det DX(x)\geq0$. The other
case is treated similarly. It turns out that for $t\in[0,1]$
sufficiently close to 1
\begin{equation}
                                             \label{9.15.3}
\int_{ D}
\det(tI+(1-t) DX(x))\,dx=\text{\rm Vol\,}D.
\end{equation}
To prove this observe that for such $t$ the Jacobian of
the mapping $X_{t}(x):=tx+(1-t)X(x)$ is 
positive on $D$ and, therefore, the image $D_{t}$ of $D$
under $X_{t}$ is a domain. For $t$ close to one also
$D_{t}\cap D\ne\emptyset$ and the mapping $X_{t}$ is invertible
(since $X(x)$ is Lipschitz continuous in $\bar D$). 
Take such a $t$.

Notice that 
if $y_{0}\in\partial D_{t}$, then there exist  $y_{n}\to y_{0}$,
$y_{n}\in D_{t}$. Then  there exist $x_{n}\in D$ such that
$y_{n}=X_{t}(x_{n})$ and for any convergent subsequence
of $x_{n}$ its limit, say $x_{0}$ is not in $D$, because
$y_{0}=X_{t}(x_{0})\not\in D_{t}$. Hence $x_{0}\in\partial D$,
$y_{0}=x_{0}$ and $\partial D_{t}\subset \partial D$. 

Similarly, if $x_{0}\in\partial D$, then there exist 
 $x_{n}\to x_{0}$, $x_{n}\in D $. Then $y_{n}:=X_{t}(y_{n})
\in D_{t}$ and $y_{n}\to y_{0}=X_{t}(x_{0})=x_{0}$. If
$y_{0}\in D_{t}$ then there is a  $z\in D$ such that
$y_{0}=X_{t}(z)=X_{t}(x_{0})$, which is impossible
since $\partial D\ni x_{0}\ne z$ and $X_{t}$ is a one-to-one
mapping in $\bar D$. Hence, $x_{0}=y_{0}\in \partial D_{t}$,
$\partial D\subset \partial D_{t}$, and $\partial D_{t}=
\partial D$.

This fact combined with the fact that
 $D$ is connected and $D_{t}\cap D\ne\emptyset$ easily implies
that $D_{t}=D$ for $t$ close to one.
 Now \eqref{9.15.3} follows. Once
\eqref{9.15.3} is true for $t$ close to 1 it is true
for all $t\in\bR$ since the left-hand side is a polynomial
with respect to $t$. By plugging in $t=1$ we get that
\begin{equation}
                                             \label{9.15.2}
\text{\rm Vol\,}D=  \int_{ D}
 \det DX(x) \,dx.
\end{equation}

Now assume that there are points $x_{0},y_{0}\in D$
such that $x_{0}\ne y_{0}$ and $z_{0}:=X(x_{0})=X(y_{0})$. 
Then there exists a (small) ball $B$ centered at $x_{0}$ which is mapped
to an open set containing $z_{0}$, such that this set is also
covered by an image of a neighborhood of $y_{0}$.
It follows that the image of $D\setminus B$ under the mapping $X$
is still $D$. Then \eqref{9.15.1} applied to $D\setminus B$
in place of $D$ shows that $\text{\rm Vol\,}D$ is less than
or equal to
the integral of $\det DX$ over $D\setminus B$ which is 
strictly less than the right hand side of \eqref{9.15.2}
since $\det DX\not\equiv0$ in $B$, because the said
neighborhood of $z_{0}$ has nonzero volume.
This is a desired contradiction and the lemma is proved.

We now know that, for each $(\omega,t)\in\Omega\times[0,T]$, 
the mapping $x\to X_{t}(x)$ is one-to-one and onto
and there exists the inverse mapping $X^{-1}_{t}(x)$, which
is infinitely differentiable in $x$ by the implicit
function theorem. In addition, from formulas for derivatives
of $X^{-1}_{t}(x)$ we conclude that these derivatives
are continuous and bounded as functions of $(t,x)$
for each $\omega$. 

Next, define the operations  ``hat'' and ``check''
which transform  any
  function
$\phi_{t}(x)$ into
\begin{equation}
                                                  \label{8.2.3}
\hat{\phi}_{t}(x):=\phi_{t}( X_{t}(x)),\quad \check{\phi}
=\phi_{t}( X^{-1}_{t}(x)).
\end{equation}
Also define $\rho_{t}(x)$ from the equation
$$
\rho_{t}(X_{t}(y))\det DX_{t}(y)=1
$$
and observe that by the change of variables formula
\begin{equation}
                                                  \label{6.30.1}
\int_{\bR^{d}}F(X_{t}(y))\phi(y)\,dy=
\int_{\bR^{d}}F(x)\check{\phi}_{t}(x)\rho_{t}(x)
\,dx,
\end{equation}
 whenever at least one side of the equation
makes sense.

We are going to make change of variables $x\to X_{t}(x)$ in 
\eqref{6.28.1} and therefore we need to understand
how the equation transforms under this change.
Define the mapping   ``bar'' which transforms    
any $\bR^{d}$-valued function $\sigma_{t}(x)$ into
\begin{equation}
                                             \label{7.23.1}
\bar{\sigma}_{t}(x)=Y_{t}(x)\hat{\sigma}_{t}(x)
=Y_{t}(x)\sigma_{t}( X_{t}(x)),
\end{equation}
where
$$
Y=(DX)^{-1}.
$$

Observe that for real-valued functions
$$
D_{j}\hat{\phi}_{t}(x)=D_{j}[\phi_{t}( X_{t}(x))]
=\widehat{D_{i}\phi}_{t}(x)D_{j}X^{i}_{t}(x),
$$
that is
$$
D\hat{\phi}=\widehat{D\phi}DX,\quad
\widehat{D\phi}=D\hat{\phi}Y.
$$

It follows that for $k=0,1,...,d_{1}+d_{2}$
\begin{equation}
                                                      \label{7.3.4}
\widehat{L_{\sigma^{k}}u}=D\hat{u}\bar{\sigma}^{k}
=L_{\bar{\sigma}^{k}}\hat{u}.
\end{equation}

One more standard fact is the following.
\begin{lemma}
                               \label{lemma 6.20.1}
For any smooth $\bR^{d}$-valued functions $\alpha$ and $\beta$
on $\bR^{d}$, for all values of arguments
we have
\begin{equation}
                                                      \label{6.30.2}
\overline{[\alpha,\beta]}=[\bar{\alpha},\bar{\beta}].
\end{equation}
\end{lemma}

Proof. Dropping the obvious values of arguments
 we have that by definition the right-hand side of \eqref{6.30.2}
equals
$$
D\bar{\beta}\bar{\alpha}-D\bar{\alpha} \bar{\beta}
=Y[\widehat{D\beta}DX\bar{\alpha}-
\widehat{D\alpha}DX\bar{\beta}]
+D_{i}Y \bar{\alpha}^{i} \hat{\beta}-
D_{j}Y \bar{\beta}^{j}\hat{\alpha}.
$$
Furthermore, since $YDX=I$,
$$
 D_{i}Y \bar{\alpha}^{i}DX+YDD_{i}X\bar a^{i}=0,
\quad D_{i}Y \bar{\alpha}^{i}=-YDD_{i}X\bar a^{i}Y,
$$
$$
 D_{i}Y \bar{\alpha}^{i} \hat{\beta}=-YDD_{i}X\bar \alpha^{i}\bar\beta
=-Y D_{ij}X\bar \alpha^{i}\bar\beta^{j}=
D_{j}Y \bar{\beta}^{j}\hat{\alpha}.
$$
This and the facts that $DX\bar{\alpha}=\hat{\alpha}$
and $DX\bar{\beta}=\hat{\beta}$ prove the lemma.

\mysection{It\^o-Wentzell formula}
                                              \label{section 6.30.4}
Here we suppose that Assumption \ref{assumption 6.28.1} (i)  
is satisfied with $R_{0}=\infty$ and define
 $C_{T}
=C_{T,\infty}$.
In this section we show 
what happens with the stochastic differential 
of a $\cD$-valued process under a random change of variables.
 
We make the following assumption which is justified in
the situation of Section \ref{section 6.30.1}
but certainly not  justified in a much more general
setting  in \cite{Ku90}.
\begin{assumption}
                                             \label{assumption 6.24.2}
There exists
a function $X_{t}(x)$ on $\Omega\times[0,T]\times\bR^{d}$ 
which has properties (i)-(iv) listed in Section \ref{section 6.30.1}
and such that, for any $(\omega,t)$, $\det DX_{t}(x)>0$
 for any $x$,
  the mapping
$x\to X_{t}(x)$ is one-to-one and onto, so that there exists an inverse
mapping $X^{-1}_{t}(x)$, and for any $R\in(0,\infty)$
$$
\sup_{\omega}\sup_{t\in[0,T]}\sup_{|x|\leq R}
|X_{t}(x)|<\infty.
$$

\end{assumption}

We start by discussing Definition \ref{def 10.25.1}
(recall that $R_{0}=\infty$).

\begin{remark}
                                       \label{remark 6.25.1}
Since $\|\cdot\|_{H^{n}_{2}}\leq \|\cdot\|_{H^{m}_{2}}$ for $n\leq m$
one can always assume that \eqref{11.16.2} holds with any $n\leq m$.
Also note that, as is well known and
easily derived by using the Fourier transform,
 for any $r\in\{0,1,...\}$ there is a constant $N$ depending only on $r$
and $d$ such that for any $\phi\in H^{2r}_{2}$
$$
\|\phi\|_{H^{2r}_{2}}:=\|(1-\Delta)^{r}\phi\|_{\cL_{2}}
\leq N\sum_{|\alpha|\leq2r}\|D^{\alpha}\phi
\|_{\cL_{2}},\quad
\sum_{|\alpha|\leq2r}\|D^{\alpha}\phi
\|_{\cL_{2}}\leq N\|\phi\|_{H^{2r}_{2}}.
$$

\end{remark}

\begin{remark}
                                                  \label{remark 6.25.2}
 Let $u\in\frD^{-\infty}_{p}(C_{T})$ and
let $\frM$ be a set of
$\cF\otimes\cB(0,T)\otimes\cB(\bR^{d})$-
measurable functions $\phi _{t}=\phi _{t}(x)=\phi _{t}(\omega,x)$
on $\Omega\times(0,T)\times\bR^{d}$ 
  such that

(i) For any $\phi\in\frM$ and $\omega$ and $t$,
 $\phi _{t}\in C^{\infty}_{0}(\bR^{d})$, and
  there exists an $R_{1} \in(0,\infty)$
 such that, for any $t\in(0,T)$, $\phi\in\frM$, and $\omega$,
we have
$\phi _{t}(x)=0$ if $|x|\geq R_{1}$;

(ii) There is an $r \in\{0,1,...\}$
such that, for $\phi\in\frM$ and $\omega\in \Omega$, the $\cL_{2}$-norm of
any derivative  of $\phi _{t}(x)$
with respect to $x$ up to order $2r$ is bounded on $(0,T) $
uniformly with respect to $\phi\in\frM$.

It turns out then that
 for any $\omega$ and any $\zeta\in C^{\infty}_{0}(\bR^{d}) $
which equals one for $|x|\leq R _{1}$ there is a constant $N$
such that for all $t\in(0,T)$
$$
 \sup_{\phi\in\frM}|(u_{t},\phi_{t} )|\leq N
\|\zeta u_{t}\|_{H^{-2r}_{2}}.
$$
In particular, if  $m$  is such that
\eqref{11.16.2} holds  
and $-r
\leq m/2$, then  
$$
\int_{0}^{T}\sup_{\phi\in\frM}|(u_{t},\phi_{t} )|^{p}\,dt<\infty.
$$

Indeed, in light of Remark \ref{remark 6.25.1}
$$
 \sup_{\phi\in\frM}|(u_{t},\phi_{t} )|=
\sup_{\phi\in\frM}|(\zeta u_{t},\phi_{t} )|\leq N
\|\zeta u_{t}\|_{H^{-2r}_{2}}\sup_{\phi\in\frM}\|\phi_{t} \|_{H^{2r}_{2}}.
$$
\end{remark}

We use the notation from Section
\ref{section 6.30.1} and observe that
$\rho_{t}(x)$ is infinitely differentiable with respect to $x$
and for any $\omega$ any its derivatives are bounded
on $[0,T]\times B_{R}$ for any $R\in(0,\infty)$.
Hence,
 the following definition makes sense:
for $u_{t}\in\cD$, $\phi\in C^{\infty}_{0}(\bR^{d})$, and $t\in[0,T]$ let
\begin{equation}
                                                      \label{6.24.1}
(\hat{u}_{t},\phi):=
(u_{t} ,\check{\phi}_{t} \rho_{t}).
\end{equation}
Observe that, if $u_{t}$ is a locally summable function,
this definition coincides with the one given in \eqref{8.2.3}
due to \eqref{6.30.1}.

\begin{lemma}
                                                \label{lemma 6.24.1}
If  $u\in\cD$ and $t\in[0,T]$, then
$(u  ,\check{\phi}_{t} \rho_{t})$ is a generalized function
for each $\omega$. Furthermore, if 
$u\in\frD^{-\infty}_{p}(C_{T})$,
then $\hat{u}\in\frD^{-\infty}_{p}(C_{T})$.
\end{lemma}

Proof. To prove the first assertion observe that
if $\phi^{n}$ converge to $\phi$ as test functions,
then their supports are in the same compact set
and $\phi^{n}\to\phi$ uniformly on $\bR^{d}$
along with each derivative in $x$. From calculus we conclude
that the same is true about $\check{\phi}^{n}_{t}\rho_{t}(x)$
for each $t$ and $\omega$ and then by definition 
$(u  ,\check{\phi}^{n}_{t} \rho_{t})\to
(u  ,\check{\phi}_{t} \rho_{t})$.

To prove the second assertion, first of all 
take a $\zeta\in C^{\infty}_{0}(\bR^{d})$ with unit integral,
define $\zeta^{n}(x)=n^{d}\zeta(n x)$ and let $u^{n}_{t}=
u_{t}*\zeta^{n}$. One knows that $u^{n}_{t}(x)$
is an infinitely differentiable function of $x$
 for each $n$, $t$, and
$\omega$
and $u^{n}_{t}\to u_{t}$ as $n\to\infty$ in the sense of 
generalized functions
 for each  $t$  and
$\omega$. In particular,
$$
(u_{t} ,\check{\phi}_{t} \rho_{t})=\lim_{n\to\infty}
(u^{n}_{t} ,\check{\phi}_{t} \rho_{t})
=\lim_{n\to\infty}\int_{\bR^{d}}u^{n}_{t} (x)
\check{\phi}_{t}(x) \rho_{t}(x)\,dx.
$$
This formula and the fact, that  for each $x$  the function
$u^{n}_{t}(x)$, continuous in $x$, is
predictable  by definition, 
show that $\hat{u}_{t}$ possesses the measurability
properties required in Definition 
\ref{def 10.25.1}.

Next, take an open ball $B\subset\bR^{d}$
and take
$\phi\in C^{\infty}_{0}(B)$.
Observe that by assumption there is an $R\in(0,\infty)$
such that $ X_{t}(x)\in B_{R}$
  for all $t\in[0,T]$, $x\in B$, and $\omega$. Take an 
$r\in\{0,1,...\}$ such that $-r\leq m/2$,
where $m$ is taken  from Definition 
\ref{def 10.25.1} corresponding to the ball $B_{2R}$, and let
$$
\frM=\{\psi\in C^{\infty}_{0}(\bR^{d}):
\|\psi\|_{H^{2r}_{2}}=1\}.
$$

Since the inequality
$\check{\phi}_{t}(x)\ne0$ implies that $X^{-1}_{t}(x)\in B$,
that is $x\in X_{t}(B) $ and $x\in B_{R}$, the supports
of   $\check{\phi}_{t}\check{\psi}_{t}$ lie in
$\bar{B}_{R}$ for all $t\in(0,T)$ and $\psi\in\frM$. It follows by
Remark \ref{remark 6.25.2} that
$$
\|\hat{u}_{t}\phi\|_{H^{-2r}_{2}}
=\sup_{\psi\in\frM}|(\hat{u}_{t},\phi\psi)|
=\sup_{\psi\in\frM}|(u_{t},\check{\phi}_{t}\check{\psi}_{t})|
\leq N\|\zeta u_{t}\|_{H^{-2r}_{2}},
$$
where $N$ is independent of $t$, and $\zeta$
is any function of class $C^{\infty}_{0}(B_{2R})$
which equals one on $B_{R}$. 
This obviously shows that $\hat{u}_{t}$ satisfies
the condition related to 
\eqref{11.16.2} if $u_{t}$ does, and the lemma is proved. 

Here is the version of It\^o-Wentzell formula we need.

\begin{theorem}
                                   \label{theorem 11.16.5} 
$ f \in\mathfrak{D}^ {-\infty}_{1}(C_{T})$, 
$u,g^{k}\in
\mathfrak{D}^ {-\infty}_{2}(C_{T})$, $k=1,...,d_{1}$,
and assume that \eqref{11.16.3} holds
 (in the sense of distributions). Then
$$
d\hat{u} _{t}= [\hat{f}_{t}+\hat{a}^{ij}_{t}
\widehat{D_{ij}u}_{t}-\hat{b}^{i}_{t}\widehat{D_{i}u}_{t}
- \widehat{D_{i}g_{t}^{k}}  \hat{\sigma}_{t}^{ik}   ]\,dt
$$
\begin{equation}
                                                  \label{4.4.5}
+
[\hat{g}^{k}_{t}-\widehat{D_{i} u}_{t}\hat\sigma^{ik}_{t} ]\,dw_{t}^{k},
\quad t\leq T
\end{equation}
(in the sense of distributions), where 
$$
a^{ij}_{t}=(1/2)\sum_{k=1}^{d_{1}}
\sigma^{ik}_{t}\sigma^{jk}_{t}.
$$
\end{theorem}

Proof. Take an   $\eta\in C^{\infty}_{0}(\bR^{d})$
and fix a $y\in\bR^{d}$. Then by Theorem 1.1 of \cite{Kr11}
the equation
$$
d(u_{t},\eta(\cdot+X_{t}(y)))
=\big([g^{k}_{t}-D_{i}u_{t}\sigma^{ik}_{t}(X_{t}(y))],
\eta(\cdot+X_{t}(y))\big)\,dw^{k}_{t}
$$
$$
\big(\big[ f_{t}+a^{ij}_{t}(X_{t}(y))D_{ij}u_{t} 
-b^{i }_{t}(X_{t}(y))D_{i}u_{t}
$$
\begin{equation}
                                                         \label{6.23.1}
- D_{i}g_{t}^{k}\sigma^{ik}_{t}(X_{t}(y)) \big],
\eta(\cdot+X_{t}(y))\big)\,dt
\end{equation}
holds, after being integrated from $0$ to $t$,
with probability one for all $t\in[0,T]$.

Then we take a $\phi\in C^{\infty}_{0}(\bR^{d})$,
 multiply both parts of
\eqref{6.23.1} by $\phi(y)$,
 and apply usual and stochastic Fubini's theorems
(see, for instance, \cite{Kr11}). Owing to the fact that for each $\omega$ and $R>0$ the set
  $\{X_{t}(y):t\in[0,T],|y|\leq R\}$ is bounded,
in order to be able to apply Fubini's theorems it suffices to show that
for any $R>0$ (a.s.)
$$
\int_{0}^{T}\sup_{|x|\leq R}
(|G_{t}(x)|+\sum_{k}|H^{k}_{t}(x)|^{2} )\,dt<\infty,
$$
where
$$
G_{t}(x)=\big(\big[ f_{t}+a^{ij}_{t}(x)D_{ij}u_{t} 
-b^{i }_{t}(x)D_{i}u_{t}- 
D_{i}g_{t}^{k}\sigma^{ik}_{t}(x) \big],
\eta(\cdot+x)\big),
$$
$$
H^{k}_{t}(x)=\big([g^{k}_{t}-D_{i}u_{t}\sigma^{ik}_{t}(x)],
\eta(\cdot+x)\big).
$$
The fact that all terms entering $G$ and $H$ apart from one
 admit needed estimates easily follows from Remark
\ref{remark 6.25.2}. The remaining one is
$$
\sum_{k}\int_{0}^{T}\sup_{|x|\leq R}
\big( D_{i}u_{t}\sigma^{ik}_{t}(x) ,
\eta(\cdot+x)\big)^{2}\,dt\leq
N \sup_{t\leq T,|x|\leq R}
\big( D u_{t}  ,
\eta(\cdot+x)\big)^{2}, 
$$
where $N<\infty$ and the last supremum is finite (a.s.)
by Lemma 4.1 of \cite{Kr11}.

Thus, we are in the position to apply Fubini's theorems.
We also use \eqref{6.30.1}.
Then we obtain
$$
d\int_{\bR^{d}}(u_{t},\eta(\cdot+x))\check{\phi}_{t}(x)\rho_{t}(x)
\,dx
$$
$$
= \int_{\bR^{d}}\big([g^{k}_{t}-
D_{i}u_{t}\sigma^{ik}_{t}(x),\eta(\cdot+x)
\big)\check{\phi}_{t}(x)\rho_{t}(x)
\,dx\,dw^{k}_{t}
$$
$$
+\int_{\bR^{d}}\big(\big[ f_{t}+a^{ij}_{t}(x)D_{ij}u_{t} 
-b^{i }_{t}(x)D_{i}u_{t}
$$
\begin{equation}
                                                         \label{6.23.2}
- D_{i}g^{k}_{t} \sigma^{ik}_{t}(x) \big],
\eta(\cdot+x)\big)
\check{\phi}_{t}(x)\rho_{t}(x)\,dx\,dt.
\end{equation}

We substitute here $\eta^{n}$ in place of $\eta$,
where $\eta^{n} $ tend to the delta-function as $n\to\infty$ in the sense
of distributions. Then we use the simple fact (having 
very little to do with
Fubini's theorem) that
$$
\int_{\bR^{d}}(u_{t},\eta^{n}(\cdot+x))\check{\phi}_{t}(x)\rho_{t}(x)
\,dx=\big(u_{t},\int_{\bR^{d}} \eta^{n}(\cdot+x))\check{\phi}_{t}(x)\rho_{t}(x)
\,dx\big),
$$
where, for each $\omega$, the test functions
$$
\int_{\bR^{d}} \eta_{n}(y+x))\check{\phi}_{t}(x)\rho_{t}(x)
\,dx
$$
as functions of $y$ vanish outside the same ball and converge to
$\check{\phi}_{t}(y)\rho_{t}(y)$ uniformly on $\bR^{d}$ along with
each derivative.
Similar statements are true about other terms entering \eqref{6.23.2},
for instance,
$$
\int_{\bR^{d}}\big( 
D_{i}u_{t}\sigma^{ik}_{t}(x), \eta^{n}(\cdot+x)
\big)\check{\phi}_{t}(x)\rho_{t}(x)
\,dx
$$
$$
=\big( 
D_{i}u_{t},\int_{\bR^{d}}\eta^{n}(\cdot+x)\sigma^{ik}_{t}(x)
\check{\phi}_{t}(x)\rho_{t}(x)
\,dx\big).
$$

We want to use the dominated convergence theorem to
pass to the limit in \eqref{6.23.2} with $\eta^{n}$
in place of $\eta$. Notice that the supports
of $\check{\phi}_{t}(y)$ and
\begin{equation}
                                               \label{11.2.1}
\int_{\bR^{d}}\eta^{n}(y+x)\sigma^{ik}_{t}(x)
\check{\phi}_{t}(x)\rho_{t}(x)
\,dx
\end{equation}
lie in the same ball for all $\omega,t ,n$.
By Remark
\ref{remark 6.25.2}
for any $\omega$
\begin{equation}
                                                         \label{6.26.1}
|\big( 
u_{t},D_{i}\int_{\bR^{d}}\eta^{n}(\cdot+x)\sigma^{ik}_{t}(x)
\check{\phi}_{t}(x)\rho_{t}(x)
\,dx\big)|^{2}\leq N\|\zeta u_{t}\|^{2}_{H^{-2r}_{2}}
\end{equation}
with $N$ independent of $t$ and $n$ if 
$\zeta\in C^{\infty}_{0}(\bR^{d})$
equals one on the supports of \eqref{11.2.1}.

 Since 
$u\in \frD^{-\infty}_{2}(C_{T})$,
the right-hand side of \eqref{6.26.1} has 
finite integral over $[0,T]$ if   $r$ is chosen appropriately,
and this allows us to pass to the limit in the stochastic
integral containing $D_{i}u_{t}$. Similarly one deals with
other integrals and concludes that
 $$
d \big(u_{t}, \check{\phi}_{t} \rho_{t} \big)=  \big([g^{k}_{t}-
D_{i}u_{t}\sigma^{ik}_{t}], 
\check{\phi}_{t} \rho_{t} 
\big)\,dw^{k}_{t}
$$
\begin{equation}
                                                         \label{6.23.3}
+ \big(\big[ f_{t}+a^{ij}_{t} D_{ij}u_{t} 
-b^{i}_{t} D_{i}u_{t}
- D_{i}g^{k}_{t}\sigma^{ik}_{t}  \big],
\check{\phi}_{t} \rho_{t}  \big)\,dt.
\end{equation}
This yields \eqref{4.4.5} by definition and the theorem is proved.

\begin{corollary}
                          \label{corollary 8.2.1}
Assume that $u$ satisfies \eqref{6.28.1}  in $C_{T}$. Then
$$
 d\hat{u}_{t} = \big[\sum_{k=1}^{d_{2}}
 L_{\bar \sigma_{t}^{d_{1}+k}}^{2} \hat{u}_{t}
+\hat{c}_{t}\hat{u}_{t}+\hat{f}_{t} -
\widehat{D_{i}g_{t}^{k}}  \hat{\sigma}_{t}^{ik} \big]\,dt
+[\hat{u}_{t}\hat{\nu}^{k}_{t}+\hat{g}^{k}_{t}]\,dw^{k}_{t}.
$$
\end{corollary}  

Indeed, as is easy to see
$$
a^{ij}_{t}D_{ij}u_{t}-b^{i}_{t}D_{i}u_{t}=(1/2)\sum_{k=1}^{d_{1}}
L_{\sigma^{k}_{t}}^{2}u_{t}-L_{\sigma^{0}_{t}} u_{t},
$$
$$
-D_{i}L_{\sigma^{k}_{t}}u\sigma^{ik}=-\sum_{k=1}^{d_{1}}
L_{\sigma^{k}_{t}}^{2}u_{t},\quad
L_{\sigma^{k}_{t}}u_{t}-D_{i}u_{t}\sigma^{ik}_{t}=0,
$$
and thanks to \eqref{7.3.4}
$$
\widehat{L_{\sigma_{t}^{d_{1}+k}}^{2} u_{t}}=
L_{\bar \sigma_{t}^{d_{1}+k}}^{2} \hat{u}_{t}.
$$

\mysection{Proof of Theorem \protect\ref{theorem 6.30.1}}
                                \label{section 6.30.2}

\begin{remark}
                                        \label{remark 7.22.1}
While proving Theorem \ref{theorem 6.30.1} we may assume
that $u_{S }=0$. Indeed, take an $s _{0}\in(S ,T)$
and take an infinitely differentiable function $\chi_{t}$, $t\geq0$,
such that $\chi_{t}=0$ on $[0,S ]$ and $\chi_{t}=1$
for $t\geq s _{0}$. Then the function $\chi_{t}u_{t}$
satisfies an easily derived equation and equals zero
at $t=S $. Furthermore, $f_{t}$ and $g^{k}_{t}$
remain unchanged for $t\geq s _{0}$ under this change
of $u$ and hence if the theorem is true when $u_{S }=0$,
then in the general case its assertions are true if
we replace in them $S $ with $s _{0}$. Due to the arbitrariness
of $s _{0}$, then the theorem is true as it is stated.

Our next observation is that, while proving 
Theorem \ref{theorem 6.30.1}, we may assume that
$S =0$. Indeed, if is not,
we can always make an appropriate shift of the origin
of the time axis.
\end{remark}

In light of Remark \ref{remark 7.22.1} everywhere
below we assume that $S =0$ and $u_{0}=0$.
The rest of the proof we split into a few steps.

{\em Step 1}. First suppose that, for $k=1,...,d_{1}$,
  $g^{k}_{t}(x)=0$ if $|x|< R_{0}$ and $\nu^{k}_{t}\equiv0$.
 Also
suppose that,
  for any $k=0,1,...,d_{1}$,  we have
$ \sigma^{ k}_{t}(x)=0$ if $|x|\geq 2R_{0}$
and $f_{t}(x)=u_{t}(x)=0$, 
 if $|x|>R_{0}-\varepsilon$, where the constant 
$\varepsilon>0$. Then equation  \eqref{6.28.1} holds on 
$ C_{T} $ in the sense of
Definition \ref{def 10.25.3} with
$\nu^{k}_{t}\equiv g^{k}_{t}\equiv0$  for $k=1,...,d_{1}$.

By Corollary \ref{corollary 8.2.1}
\begin{equation}
                                                 \label{7.3.5}
 (\hat{u}_{t},\phi)= \int_{0}^{t}\big(\sum_{k=1}^{d_{2}}
 L_{\bar \sigma_{s}^{d_{1}+k}}^{2} \hat{u}_{s}
+\hat{c}_{s}\hat{u}_{s}+\hat{f}_{s},\phi\big)\,ds
\end{equation}
and this holds for any $\phi\in C^{\infty}_{0}(\bR^{d})$ 
with probability one
for all $t\in[0,T]$. Let $\Phi$ be a countable subset of
$C^{\infty}_{0}(\bR^{d})$ which is everywhere dense in  $H^{n}_{2}$
for any $n\in\bR^{d}$. Then there exists a set $\Omega'$
of full  probability such that for any $\omega\in\Omega'$
and any $\phi\in\Phi$ equation \eqref{7.3.5} holds
for all $t\in[0,T]$.
 By setting $u$ and $f$ to be zero
if necessary
for $\omega\not\in\Omega'$ we may assume that
equation  \eqref{7.3.5} 
 holds for any $\phi\in\Phi$,
 $t\in[0,T]$, and $\omega$. Furthermore, observe that
by assumption $u_{t}(x)=0$ if $|x|\geq R_{0}-\varepsilon$.
 Hence, \eqref{11.16.2}
holds with $\phi\equiv1$ with probability one for an 
appropriate $m$ (depending on $\omega$).
 By redefining,
if necessary, $u$ one more time we may assume that for any $\omega$
there exists an integer $r$ such that
\begin{equation}
                                                \label{7.3.6}
\int_{0}^{T}\|u_{t}\|_{H^{-2r}_{2}}^{2}\,dt<\infty,\quad
\int_{0}^{T}\|\hat{u}_{t}\|_{H^{-2r}_{2}}^{2}\,dt<\infty.
\end{equation}
Having this and similar relations for $f$ and
remembering that $\Phi$ is dense in $H^{2r}_{2}$ we easily conclude that
\eqref{7.3.5} holds for any $\phi\in C^{\infty}_{0}(\bR^{d})$,
 $t\in[0,T]$, and $\omega$.

Next argument is conducted for a fixed $\omega\in\Omega_{0}$. 
Introduce
$$
\hat{G}=\{(t,x):t\in(0,T),x\in X^{-1}_{t}(B_{R_{0}})\}.
$$
Since $X_{t}(x)$ is a diffeomorphism continuous with respect to $t$,
$\hat{G}$ is a domain.  Furthermore, it follows from the assumptions
of the theorem that for any $\zeta\in C^{\infty}_{0}(\hat{G})$ and any
$n=1,2,...$, we have
$$
\int_{0}^{T}\|\hat{f}_{t}\zeta\|^{2}_{H^{n}_{2}}\,dt<\infty.
$$

Next let
$\bar{\,\bL}_{0}=\{\bar \sigma^{d_{1}+1},...,\bar \sigma^{{d_{1}+d_{2}}}\}$,  
$$
\bar{\,\bL} _{n+1}=\bar{\,\bL}_{n}\cup\{[ 
\bar\sigma^{d_{1}+k},M]:k=1,...,d_{2},
M\in\bar{\,\bL}_{n}\},\quad n\geq0.
$$ 
Note that by Lemma \ref{lemma 6.20.1}, if $\sigma\in\bL_{n}$,
then $\bar\sigma\in\bar\bL_{n}$.

Now, take   $\zeta \in C^{\infty}_{0}(\hat G)$ and $\zeta_{1}
\in C^{\infty}_{0}( G)$ so that 
$$
\zeta_{1}=1 \quad\text{on}\quad \supp \check\zeta .
$$
 By Assumption \ref{assumption 1.5.1}
there exists an $n\in\{0,1,...\}$ such that for any  
 $i=1,2,...,d$ there exist 
$r\in\{0,1,...\}$ and $\sigma^{(i1)},...,\sigma^{(ir)}\in \bL_{n}$
and real-valued functions $\gamma^{(i1)},...,\gamma^{(ir)}$
of class $BC^{\infty}_{b}$ such that
$$
\zeta_{1} e_{i}=\gamma^{(i1)}\sigma^{(i1)}+...+\gamma^{(ir)}\sigma^{(ir)}.
$$
Obviously one may  assume that $r$ is common for all $i=1,2,...,d$.  
It follows that
$$
\hat \zeta_{1} Ye_{i}=\hat\gamma^{(i1)}
\bar\sigma^{(i1)}+...+\hat \gamma^{(ir)}\bar\sigma^{(ir)},
$$
which after being multiplied by $\zeta$ yields
$$
 \zeta  Ye_{i}=\zeta \hat\gamma^{(i1)}
\bar\sigma^{(i1)}+...+\zeta \hat \gamma^{(ir)}\bar\sigma^{(ir)},
$$
Observe that for $\xi\in\bR^{d}$ and $\lambda=DX\xi$ 
we have $Ye_{i}\lambda^{i}=\xi$, so that
$$
 \zeta \xi=\lambda^{i}\hat\gamma^{(i1)}
\bar\sigma^{(i1)}+...+\lambda^{i}\hat \gamma^{(ir)}\bar\sigma^{(ir)}.
$$
Hence, for any $\xi\in\bR^{d}$ and $\zeta\in C^{\infty}_{0}(\hat G)$,
$\zeta\xi$ is represented as a linear combination
of elements of $\bar{\,\bL}_{n}$ with coefficients of 
class $BC^{\infty}_{b}$

We checked the assumptions of
  Theorem 2.7 of \cite{Horm} and by that theorem conclude
that $\hat{u}_{t}(x)$
is infinitely differentiable with respect to $x$ for
$(t,x)\in \hat{G}$, each of its 
derivative is a continuous function
in $\hat{G}$, and an estimate similar to
\eqref{7.22.2} holds. Changing back the coordinates
we get the 
first assertion of our theorem
and, in addition, the fact that
in \eqref{7.16.1} the right-hand side can be taken  to
be $d$ in place of $d+1$.

{\em Step 2}. We keep the assumption  of Step 1 that,
  for $k=1,...,d_{1}$,
 $g^{k}_{t}(x)=0$ if $|x|< R_{0}$ and  $\nu^{k}_{t}\equiv0$.
We will cut-off $u_{t}$ for $x$ near the boundary of
$B_{R_{0}}$, so that the new function will satisfy
an equation in $C_{T}$, to which we can then
apply the It\^o-Wentzell formula.
The only difficulty which appears after that
is that we will get a new $g^{k}_{t}$ which is not 
vanishing in $C_{T,R_{0}}$. Partial help comes 
from the fact that if we cut-off close to the boundary,
then the new $g^{k}_{t}$ will be not vanishing
only near the boundary. Due to this fact the transformations
made in Step 1 will not lead exactly to
a deterministic equation like \eqref{7.3.5}
with random coefficients but to an equation
containing the  stochastic
integral of $\hat{g}^{k}_{t}\,dw_{t}^{k}$. This integral
can be, so to speak, {\em locally in time}
 neglected near the lateral boundary of a domain like
$\hat{G}$. This yields a deterministic situation
where we apply Theorem 2.7 of \cite{Horm}.

Take a sequence $\zeta^{n}\in C^{\infty}_{0}(B_{R_{0}})$ such that
$\zeta^{n}=1$ on $B_{R_{0}-1/n}$ and $\zeta^{n}=0$ on $B_{R_{0}-1/(n+1)}$
and define
$u^{n}_{t}=u_{t}\zeta^{n}$. Then as is easy to see
\begin{equation}
                                                    \label{7.13.1}
du^{n}_{t}=(L_{t}u^{n}_{t}
+c_{t}u^{n}_{t}+f^{n}_{t})\,dt
+(L_{\sigma^{ k}_{t}}u^{n}_{t}+g^{nk}_{t} )\,dw^{k}_{t} 
\end{equation}
in $C_{T}$, where 
$$
f^{n}_{t}=f_{t}\zeta^{n}-u _{t}L_{t}\zeta^{n}
- (L_{\sigma^{ k}_{t}}u _{t})
L_{\sigma^{ k}_{t}} \zeta^{n},\quad g^{nk}_{t}=
-u_{t}L_{\sigma^{ k}_{t}}\zeta^{n}.
$$

Also take a $\zeta\in C^{\infty}_{0}(\bR^{d})$ such that
$\zeta=1$ on $B_{R_{0}}$ and $\zeta=0$ outside $B_{2R_{0}}$.
Obviously, in \eqref{7.13.1} one can replace the operator
$L_{t}$  
  with the one  denoted by 
$\tilde L _{t}$   and constructed on the basis
of $\tilde \sigma^{ k}_{ t}:=\zeta \sigma^{ k}_{t}$.
Thus,
\begin{equation}
                                                    \label{7.13.2}
du^{n}_{t}=(\tilde L _{t}u^{n}_{t}
+c_{t}u^{n}_{t}+f^{n}_{t})\,dt
+( L _{\tilde\sigma^{ k}_{t}}u^{n}_{t}+g^{nk}_{t} )\,dw^{k}_{t},
\end{equation}

Next we change the coordinates by defining $X _{ t}(x)$
as a unique solution of
\begin{equation}
                                                 \label{7.23.2}
x_{t}=-\int_{0}^{t}\tilde \sigma^{k}_{s}( x_{s})\,dw^{k}_{s}
-\int_{0}^{t}\tilde b_{t}(x_{s})\,ds,
\end{equation}
where
$$
\tilde b_{t}(x)=\tilde \sigma^{0}_{t}(x)-(1/2)
\sum_{k=1}^{d_{1}}D\tilde \sigma^{k}
_{t}(x)\tilde \sigma^{k}_{t}(x).
$$
We also recall that $u\in \frD^{-\infty}_{2}
(C_{T,R_{0}})$ so that the stochastic
integral 
$$
m^{n}_{t}:=\int_{0}^{t}\hat {u}_{t}
\widehat{L^{n}_{\sigma^{ k}_{s}}\zeta^{n}}\,dw^{k}_{s}
$$
is well-defined as a stochastic integral  of a Hilbert-space
valued function and is continuous with respect to $t$
for all $\omega$. Then similarly to
\eqref{7.3.5} we come to the conclusion that
for any $\phi\in C^{\infty}_{0}(\bR^{d})$ with probability
one 
\begin{equation}
                                                 \label{7.13.02}
 (\hat{u}^{n}_{t},\phi)=(\hat{u}_{0},\phi)+\int_{0}^{t}\big(\sum_{k=1}^{d_{2}}
 L_{\bar \sigma_{ s}^{d_{1}+k}}^{2} \hat{u}^{n}_{s}
+\hat{c}_{s}\hat{u}^{n}_{s}+\hat{f}^{n}_{s},\phi\big)\,ds
-(\phi,m^{n}_{t})
\end{equation}
for all $t\in[0,T]$,
where $\bar \sigma_{ s}^{ k}$ are constructed
from $\tilde\sigma_{ s}^{ k}$ as in \eqref{7.23.1}
starting with equation \eqref{7.23.2} instead of
\eqref{6.28.4}.

After that by doing the same manipulations as below \eqref{7.3.5}
we convince  ourselves that without losing generality
we may assume that \eqref{7.13.02} holds for all
 $\phi\in C^{\infty}_{0}(\bR^{d})$,
 $t\in[0,T]$, and $\omega$. This and our result about 
$\bar\bL_{n}$ 
are the only
facts which we need from the arguments in Step 1.

Then we again argue with $\omega\in\Omega_{0}$ fixed.
Take $t_{0} \in(0,T)$ and $y_{0}\in B_{R_{0}-2/n}$.
Then there is an $\varepsilon>0$
such that, for $x _{0}=X^{-1}_{t_{0}}(y_{0})$
we have 
$$
X_{ t}(B_{\varepsilon}(x _{0}))
\subset B_{R_{0}-1/n}
$$
 for any $t\in
(t_{0}-\varepsilon,
t_{0}+\varepsilon)$. For $\phi\in C^{\infty}_{0}(
B_{\varepsilon}(x _{0}))$ and $t\in
(t_{0}-\varepsilon,
t_{0}+\varepsilon)$ we have $(\phi,m^{n}_{t})=(\phi,m^{n}_{
t_{0}-\varepsilon})$ by Lemma \ref{lemma 7.13.1} since, for those $t$,
$ L_{\sigma^{k}_{t}}\zeta^{n}=0$ in $B_{R_{0}-1/n}$,
$\check{\phi}=0$ outside $B_{R_{0}-1/n}$,
$\check{\phi}L_{\sigma^{k}_{t}}\zeta^{n}\equiv0$, and
$$
(\phi,\hat {u}_{t}
\widehat{L_{\sigma^{ k}_{t}}\zeta^{n}})
=(u_{t},\rho_{t}\check{\phi}L_{\sigma^{ k}_{t}}\zeta^{n})=0.
$$

It follows that for $\phi\in C^{\infty}_{0}(
B_{\varepsilon}(x^{n}_{0}))$ and $t\in
(t_{0} -\varepsilon,
t_{0}+\varepsilon)$ 
$$
(\hat{u}^{n}_{t},\phi)=(\hat{u}_{t_{0}-\varepsilon},\phi)+
\int_{t_{0}-\varepsilon}^{t}\big(\sum_{k=1}^{d_{2}}
 L_{\bar \sigma_{n,s}^{d_{1}+k}}^{2} \hat{u}^{n}_{s}
+\hat{c}_{s}\hat{u}^{n}_{s}+\hat{f}_{t}\hat{\zeta}^{n}_{s},\phi\big)\,ds.
$$
As in Step 1 we conclude by Theorem 2.7 of \cite{Horm}
that $\hat{u}^{n}_{t}(x)$
is infinitely differentiable with respect to $x$ for
$(t,x)\in G _{\varepsilon}:=(t_{0} -\varepsilon,
t_{0}+\varepsilon)\times B_{\varepsilon}(x _{0})$ 
and each derivative is a continuous function
in $G _{\varepsilon}$. 
Furthermore, an estimate similar to \eqref{7.22.2}
is available for any closed cylinder 
inside $G _{\varepsilon}$. Actually,
Theorem 2.7 of \cite{Horm} is formally applicable only
if $t_{0} -\varepsilon=0$ and $\hat{u}_{t_{0}-\varepsilon}=0$.
Our explanations given in Remark \ref{remark 7.22.1}
take care of the general case.

Changing back the coordinates
we get that $u^{n}_{t}(y)$ is infinitely differentiable 
with respect to $y$ for
$y$ in a neighborhood of $y_{0}$ and $t$
in a neighborhood of $t_{0} $ and each derivative
is a continuous function of $(t,y)$ for those $(t,y)$.
Estimate \eqref{7.22.2} is also valid in any closed
cylinder lying in that neighborhood.
Since $y_{0}\in B_{R_{0}-2/n}$, the said neighborhood
of $y_{0}$ can be taken to belong to $B_{R_{0}-1/n}$,
where $u^{n}_{t}=u_{t}$ and $f^{n}_{t}=f_{t}$.
Now the assertion of the theorem follows 
owing to the arbitrariness
of $y_{0}$, which is provided by the possibility to take $n$
as large as we wish. Again as in Step 1 it suffices that condition
\eqref{7.16.1} be satisfied with $d$ in place of $d+1$.

{\em Step 3}. 
 Now we abandon the assumption of Step 2 
that,
  for $k=1,...,d_{1}$,
 $\nu^{k}_{t}\equiv0$, but still assume that 
$g^{k}_{t}(x)=0$ if $|x|< R_{0}$ for $k=1,...,d_{1}$.
Introduce the function $v_{t}(x,y)=yu_{t}(x)$ and the $d+1$-dimensional 
vectors
$$
\ddot{\sigma}^{k}_{t}(x,y)=\begin{pmatrix}\sigma^{k}_{t}(x)\\y
\nu^{k}_{t}(x)
\end{pmatrix},\quad k\leq d_{1},\quad
\dot{\sigma}^{k}_{t}(x,y)=
\begin{pmatrix}\sigma^{k}_{t}(x)\\0
\end{pmatrix} ,\quad k\leq d_{1}+d_{2},
$$
$$
 \dot\sigma ^{d_{1}+d_{2}+1}_{t}(x,y)=\begin{pmatrix}0\\1
\end{pmatrix}.
$$
Obviously, Assumption \ref{assumption 1.5.1} is satisfied
if we replace $G$, $d$, and $d_{2}$ 
with $G\times(0,1)$, $d+1$, and $d_{2}+1$,
respectively.
Also, routine computations yield that
 $v_{t}$ satisfies
$$
dv_{t}=\bigg((1/2)\sum_{k=1}^{d_{1}}[L^{2}_{\ddot\sigma^{k}_{t}}v_{t}
-\nu^{k}_{t}L_{\dot\sigma^{k}_{t}}v_{t}-
L_{\dot\sigma^{k}_{t}}(\nu^{k}_{t}v_{t})
-\nu^{k}_{t}\nu^{k}_{t}v_{t}]
$$
$$
+(1/2)\sum_{k=1}^{ d_{2}+1}L^{2}_{\dot\sigma^{d_{1}+k}_{t}}v_{t}
+L_{\dot\sigma^{0}_{t}}v_{t}+c_{t}v_{t}+yf_{t}\bigg)\,dt
+ L_{\ddot\sigma^{ k}_{t}}v_{t} \,dw^{k}_{t}.
$$
The result of Step 2 shows that $v_{t}$ is infinitely differentiable
with respect to $(x,y)$ in $B_{R_{0}}\times(0,1)$ for any $t\in (0,T)$
and
 the derivatives are continuous with respect to $(t,x,y)$.
Also the corresponding counterpart of \eqref{7.22.2} holds for $v_{t}$
under condition \eqref{7.16.1}. This obviously proves the theorem
in this particular case.

{\em Step 4}. Now we consider the general case.
Take an $R'_{0}\in(0,R_{0})$ and
  $\zeta\in C^{\infty}_{0}(B_{R_{0}})$
such that $\zeta=1$
on $B_{R'_{0}}$.
Then according to classical results, for sufficiently large
constant $K>0$, there exists
a  function $v\in\frD^{-\infty}_{2}(C_{T})$ such that $v_{0}=0$,
$$
\int_{0}^{T}\|v_{t}\|^{2}_{H^{n}_{2}}\,dt<\infty
$$
(a.s.) for any $n$,
 and
$$
dv_{t}=K\Delta v_{t}\,dt+
(L_{\zeta\sigma^{ k}_{t}}v_{t}+\zeta \nu^{k}_{t}v_{t}
 +\zeta g^{k}_{t})\,dw^{k}_{t}.
$$

Then the function $w_{t}=u_{t}-v_{t}$ satisfies an equation
which falls into the scheme of Step 3 with $R_{0}'$
in place of $R_{0}$ and a different $f$
but still satisfying the assumption of Theorem \ref{theorem 6.30.1}.

The assertion of the theorem now follows and the theorem is proved.

 \mysection{Proof of Theorem 
\protect\ref{theorem 7.22.1}}
                                 \label{section 8.5.1}

 The idea of the proof is to find a neighborhood of
 $[S ,T]\times B_{r}$ to which Theorem \ref{theorem 6.30.1}
is applicable. First, we extend $u_{t}$ beyond $T$.
 To do that we take 
$R_{1}\in(0,R_{0})$,
 $\zeta\in C^{\infty}_{0}(B_{R_{0}})$,
which equals one in   $  B_{R_{1}}$
 and
consider the 
function $v_{t}=\zeta u_{t}$ for $t\in[0,T]$. 
By Remark \ref{remark 11.1.1} there is an $m\in\bR$ such that,
 with probability one,
$v_{t}$ is a continuous $H^{m}_{2}$-valued process. 

It follows that $v_{T}\in H^{m }_{2}$ (a.s.), so that
solving the heat equation
$$
dv_{t}=\Delta v_{t}\,dt,\quad t> T, x\in\bR^{d}, 
$$
with initial data $v_{T}$, which is possible by classical results,
 allows us to extend $v_{t}$
beyond $T$ as an
$H^{m}_{2}$-valued  continuous functions of $t$. If we 
now  
accordingly define $c,f,\nu,g,\sigma $ for $t\geq
T$, then we will see that the assumptions of
Theorem \ref{theorem 6.30.1} are satisfied with $(
S ,T+1)\times B_{R_{1}}$ in place of $G$.
This proves the first assertion of Theorem
\ref{theorem 7.22.1}.

Passing to the second one 
we assume that $u_{S } $ is infinitely differentiable
in $B_{R_{0}}$ for every $\omega\in\Omega_{0}$. Then
we want to reduce the general case to the one
in which  $u_{S }=0$ 
in $B_{R_{0}}$ for $\omega\in\Omega_{0}$.
To achieve that take $R_{1}\in(r,R_{0})$ and
 $\zeta\in C^{\infty}_{0}(B_{R_{0}})$
as in the beginning of the proof
and solve the equation
\begin{equation} 
                                         \label{8.5.2}
dv_{t}=\big[\Delta v_{t}
+(1/2)\sum_{k=1}^{d_{1}}
L_{\zeta\sigma^{k}_{t}}^{2}v_{t}\big]\,dt+L_{\zeta\sigma^{k}_{t}}v_{t}\,dw^{k}_{t},\quad
t\in(S ,T)
,x\in \bR^{d}
\end{equation}
with initial data $v_{S }=\zeta u_{S }$. 
After making an appropriate random change of coordinates
according to Corollary \ref{corollary 8.2.1} we reduce
this SPDE to a  usual parabolic equation with random
coefficients which is uniformly
nondegenerate for any $\omega\in\Omega $
(we said more about this in  
the beginning of Section \ref{section 6.30.2}).
By classical results
there is a  solution $v_{t}$ of this new equation
  with initial data $\zeta u_{S }$,  
 which, for any $\omega\in\Omega_{0}$,
is continuous  in $[S ,T]\times \bR^{d}$
along with each its derivative of any order with 
respect to $x$. This is true because   
$\zeta u_{S }\in C^{\infty}_{0}(\bR^{d})$
for $\omega\in\Omega_{0}$. 
The same holds for equation \eqref{8.5.2}.
Furthermore, 
\begin{equation}
                                        \label{7.22.3}
\sup_{(t,x)\in[S ,T]\times \bR^{d}}
|D^{\alpha}v_{t}(x)|^{2}+
\int_{S }^{T}\|v_{t}\|_{H^{l+2}_{2}}^{2}\,dt
\leq N\|\zeta u_{S }\|^{2}_{H^{l+1}_{2}},
\end{equation}
provided that $2(l+1-|\alpha|)>d$ and $\omega\in\Omega_{0}$. 

We set $v_{t}=\zeta u_{t}$ for $t\in[0,S ]$
and then we see that   
in $[0 ,T]\times B_{R_{1}}$ the function
$u_{t}-v_{t}$ satisfies the same equation as $u_{t}$
with   $g^{k}_{t}I_{(S,T)}$ in place of $g^{k}_{t}$ and
 with  a new $f_{t}$,
 whose norms for $\omega\in\Omega_{0}$ admit and
  obvious estimates
through the norms of the old one 
and the right-hand side of \eqref{7.22.3}.
Hence, the assumptions of the present theorem
are satisfied with $ B_{R_{1}}$
in place of $B_{R_{0}}$.

By replacing $u_{t}$ and $R_{0}$
with $u_{t}-v_{t}$ and $R_{1}$,  
 we see
 that without loosing generality we may assume
that $u_{t}=f_{t}=g^{k}_{t}=0$ for $t\in[0,S ]$
on $ B_{R_{0}}$.
 In that case, we define $u_{t}=f_{t}=g^{k}_{t}
=0$, $\sigma^{k}_{t}=0$, $k=0,1,...,d_{1}+d_{2}$, for $t\in[ -1,S )$.
We also introduce new $\sigma^{k}_{t}$ for
 $k=d_{1}+d_{2}+i$, $i=1,...,d$,
by setting $\sigma^{k}_{t}=e_{i}I_{[ -1,S )}(t)$,
where the $e_{i}$'s
form the standard orthonormal basis in $\bR^{d}$. After that
we define $L_{t}$ for $t\in[ -1,S )$ according to \eqref{7.22.4},
where we replace $d_{1}+d_{2}$ with $d_{1}+d_{2}+d$ and observe that
the new $u_{t}$ now satisfies \eqref{6.28.1} in $(-1,T)
\times B_{R_{0}}$. The reader may object that $w^{k}_{t}$
are not defined for negative $t$,
but since $dw^{k}_{t}$ for negative $t$ are multiplied by zeros,
one can just take  independent Wiener processes and glue
 them to $w^{k}_{t}$ from $ -1$ to $0$. As is easy to see,
the first assumption  of the present theorem is satisfied with
$ I_{[ -1,T]}$ in place of $ I_{[S,T]}$,
 and this proves the present theorem.

\end{document}